\documentclass[number,citesort,MSNbibl,seceqn,dvips]{arxbj}

\aid{0}
\volume{17}
\issue{4}
\pubyear{2011}
\firstpage{1268}
\lastpage{1284}
\doi{10.3150/10-BEJ313}

\makeatletter
\newtheorem{tthhmm}{Theorem}[section]
\newtheorem{lem}{Lemma}[section]

\newtheorem{prop}{Proposition}[section]

\newremark{rem}{Remark} 

\newcommand{\N}{{\mathbb N}}
\newcommand{\R}{{\mathbb R}}
\newcommand{\Z}{{\mathbb Z}}

\def\bptnote#1{}

\makeatother

\begin{document}
\begin{frontmatter}

\title{Absolute regularity and ergodicity of Poisson count processes}

\runtitle{Ergodicity of Poisson count processes}

\begin{aug}
\author{\fnms{Michael H.} \snm{Neumann}\corref{}\ead[label=e1]{michael.neumann@uni-jena.de}}

\runauthor{M.H. Neumann}
\address{Friedrich-Schiller-Universit\"{a}t Jena,
Institut f\"{u}r Stochastik,
Ernst-Abbe-Platz 2, 07743 Jena,
Germany. \printead{e1}}
\end{aug}

\received{\smonth{2} \syear{2010}}
\revised{\smonth{7} \syear{2010}}

%
\begin{abstract}
We consider a class of observation-driven Poisson count processes where
the current
value of the accompanying intensity process depends on previous values of
both processes. We show under a contractive condition that the
bivariate process
has a unique stationary distribution and that a stationary version of
the count process
is absolutely regular. Moreover, since the intensities can be written
as measurable
functionals of the count variables, we conclude that the bivariate process
is ergodic.
As an important application of these results, we show how a test method
previously used in the
case of independent Poisson data can be used in the case of Poisson
count processes.
\end{abstract}

%

\end{frontmatter}

\section{Introduction}
\label{S1}

The modeling and the analysis of count data has received increasing
attention during the last decade. There are possible applications in
various fields, such as biometrics, econometrics and finance; see
Davis, Dunsmuir and Wang \cite{DDW99} and Davis and Wu \cite{DW09} for
examples. A comprehensive account of models for time series of counts
is given in Kedem and Fokianos \cite{KF02}, Chapter~4. In the majority
of cases the count variables are assumed to be Poisson distributed,
conditioned on the past and perhaps some additional regressor
variables. Models for count data consist of at least two processes: an
observable process of counts and an accompanying intensity process that
is usually not observed. Cox \cite{Cox81} and later Davis, Dunsmuir and
Wang \cite{DDW99} classified these models into parameter-driven and
observation-driven specifications. In the first case, the accompanying
intensity process evolves independently of the past history of the
observation process while, in the latter case, the values on the
intensity process do depend on past observations. The major aim of this
paper is to derive important properties such as stationarity, mixing
and ergodicity for a certain class of observation-driven processes.
Davis, Dunsmuir and Wang \cite{DDW99} mentioned that, in contrast to
parameter-driven models where these properties are inherited by the
observation process from the accompanying intensity process, there is
much less theory available in the case of observation-driven models.
Actually, ergodicity has been shown so far in a few special cases only
-- see Grunwald, Hyndman, Tedesco and Tweedie \cite{GHTT00}, Streett
\cite{Str00}, Davis, Dunsmuir and Streett \cite{DDS03}, Zheng and
Basawa \cite{ZB08} and Fokianos, Rahbek and Tj{\o}stheim \cite{FRT09}.
In these papers, the authors could use classical Markov chain theory.

In the present paper, we study a model where the observations $N_t$ are
Poisson distributed, conditioned on the past, with an intensity
$\lambda_t$ depending on one lagged value of the count process and the
intensity process; that is, $\lambda_t=f(\lambda_{t-1},N_{t-1})$, for
some function~$f$. Models of this type have been considered
before~by~Rydberg and Shephard \cite{RS00}, Streett \cite{Str00},
Davis, Dunsmuir and Streett \cite{DDS03}, Fokianos, Rahbek and
Tj{\o}stheim \cite{FRT09} and Fokianos and Tj{\o}stheim \cite{FT10}.
An important aspect is that such models allow for an autoregressive
(AR) feedback mechanism in the intensity process and it can be expected
that this leads to a parsimonious parametrization. For clarity of
exposition, we do not include additional regressor variables that are
often also incorporated in specifications of the intensity. Under a
contractive condition on~$f$, we state in Section~\ref{S2} that the
bivariate process $((N_t,\lambda_t))_{t\in\N}$ has a unique stationary
distribution. The proof of this result is based on a simple
construction, where independently started versions of the process are
coupled in such a way that they converge to each other.
Section~\ref{S3} contains the main results. For a stationary version of
the process, we prove absolute regularity ($\beta$-mixing) of the
univariate count process. Since the latter process is not Markovian, we
cannot rely on standard arguments from Markov chain theory; rather, we
use coupling arguments to derive this result. We also discuss an
example that shows that the bivariate process
$((N_t,\lambda_t))_{t\in\Z}$ and even the intensity process
$(\lambda_t)_{t\in\Z}$ are not absolutely regular in general. However,
since the intensities can be written as measurable functionals of the
count variables, we conclude from the mixing property of the count
process that the bivariate process is ergodic. In Section~\ref{S4}, we
propose a test for a particular specification of the intensity process.
We use a test statistic that has been applied before by several authors
in connection with independent Poisson random variables. Using the
ergodicity result from Section~\ref{S3}, we can show that the test
statistic is asymptotically normal. All proofs are deferred to a final
Section~\ref{SP}.

\section{Stationarity of the bivariate process}
\label{S2}

We assume that $(N_t)_{t\in\N}$ is a time series of counts, accompanied
by an intensity process $(\lambda_t)_{t\in\N}$. Denote by ${\mathcal
B}_t^{N,\lambda}=\sigma(\lambda_1,\ldots,\lambda_t, N_1,\ldots,N_t)$
the $\sigma$-field generated by the past and present values of the two
processes at time~$t$. We assume throughout that
%
\begin{equation}
\label{2.1} N_t\vert{\mathcal B}_{t-1}^{N,\lambda} \sim
\operatorname{Poisson}(\lambda_t)
\end{equation}
and
%
\begin{equation}
\label{2.2} \lambda_t = f(\lambda_{t-1},N_{t-1})
\end{equation}
for some function $f\dvtx  [0,\infty) \times\N_0\rightarrow[0,\infty)$
($\N_0=\N\cup\{0\}$). For the time being, the starting value
$\lambda_1$ may be random or non-random. It follows from the structure
of the model that ${\mathcal
B}_{t-1}^{N,\lambda}=\sigma(\lambda_1,N_1,\ldots,N_{t-1})$ and that the
bivariate process $((N_t,\lambda_t))_{t\in\N}$ forms a homogeneous
Markov chain. Throughout this paper we will assume that the
function~$f$ satisfies the following contractive condition:
%
\begin{equation}
\label{2.3} | f(\lambda,y) - f(\lambda',y') | \leq
\kappa_1 |\lambda-\lambda'| + \kappa_2 |y-y'| \qquad\forall
\lambda,\lambda'\geq0, \forall y,y'\in\N_0,
\end{equation}
where $\kappa_1$ and $\kappa_2$ are non-negative constants with
$\kappa:=\kappa_1+\kappa_2<1$. This includes as a special case a
linear specification where $\lambda_t=\theta_0+\theta_1
\lambda_{t-1}+\theta_2 N_{t-1}$ and $\theta_0,\theta_1,\theta_2$ are
non-negative constants with $\theta_1+\theta_2<1$. Rydberg and Shephard
\cite{RS00} proposed such a model for describing the number of trades
on the New York Stock Exchange in certain time intervals and called it
the $\operatorname{BIN}(1,1)$ model. Stationarity and other properties for this model
were derived by Streett \cite{Str00} and Ferland, Latour and Oraichi
\cite{FLO06}, who referred to it as the $\operatorname{INGARCH}(1,1)$ model, and
Fokianos, Rahbek and Tj{\o}stheim \cite{FRT09}. The generality of
Condition (\ref{2.3}) is chosen to include nonlinear specifications
such as the exponential AR model proposed in Fokianos, Rahbek and
Tj{\o}stheim \cite{FRT09}. In this case, the intensity function is
specified as $f(\lambda,y)=(a+c \exp(-\gamma\lambda^2))\lambda+by$,
where $a,b,c,\gamma>0$. It follows from $\frac{\partial}{\partial y}
f(\lambda,y)=b$ and $|\frac{\partial}{\partial\lambda}
f(\lambda,y)|\leq a+c$ that (\ref{2.3}) is fulfilled if $a+b+c<1$.

Note that (\ref{2.3}) implies that
%
\begin{equation}
\label{2.4} f(\lambda,y) \leq f(0,0) + \kappa_1 \lambda+
\kappa_2 y.
\end{equation}
It follows from (\ref{2.4}) that $E(\lambda_t|\lambda_{t-1})\leq
f(0,0)+\kappa\lambda_{t-1}$, which leads to
%
\begin{equation}
\label{2.5} E( N_t|\lambda_1 ) = E( \lambda_t|\lambda_1 )
\leq f(0,0) \frac{1-\kappa^{t-1}}{1-\kappa} + \kappa^{t-1}
\lambda_1.
\end{equation}
Hence, the bivariate chain $((N_t,\lambda_t))_{t\in\N}$ is bounded in
probability on average. Moreover, it follows from (\ref{2.3}) that, for
any open set $O\in2^{\N_0}\otimes{\mathcal B}$, the transition
probabilities $P((N_t,\lambda_t)\in O|(N_{t-1},\lambda_{t-1})=
\cdot)$ are a continuous, and therefore also a lower semicontinuous,
function. Hence, the Markov chain is a weak Feller chain and it follows
from Theorem~12.1.2(ii) in Meyn and Tweedie \cite{MT93} that there
exists at least one stationary distribution. Uniqueness of this
stationary distribution, however, requires more than (\ref{2.4}) and
will follow from the contractive condition (\ref{2.3}). The following
theorem summarizes this and a few other useful facts.

\begin{tthhmm}\label{T2.1}
Suppose that the bivariate chain
$((N_t,\lambda_t))_{t\in\N}$ obeys (\ref{2.1})--(\ref{2.3}). Then
\begin{longlist}[(iii)]
\item There exists a unique stationary distribution~$\pi$.
\item If $(N_1,\lambda_1)\sim\pi$, then $E\lambda_1<\infty$.
\item If $f(0,0)=0$, then $\pi(\{0,0\})=1$.
If $f(0,0)>0$, then $\pi(\{y,\lambda\})<1$ for all $y\in\N_0$,
$\lambda\in[0,\infty)$.
\end{longlist}
\end{tthhmm}

\begin{rem}
Using the contractive property (\ref{2.3}), we will show in the
proof of Theorem~\ref{T2.1} that the~$n$-step transition laws
$P((N_{n+1},\lambda_{n+1})\in\cdot|(N_1,\lambda_1)=x)$ converge to
a common limit $\pi$ not depending on the starting value~$x$, where
$\pi$ is a probability measure. This will imply that~$\pi$ is the
unique stationary distribution.

There are alternative ways to prove Theorem~\ref{T2.1}. Introducing a
sequence of independent ``innovations'' $(U_t)_{t\in\N}$ with
$U_t\sim\operatorname{Uniform}[0,1],$ we could re-express the process values as
\[
(N_{t+1},\lambda_{t+1}) = G((N_t,\lambda
_t),U_{t+1}) :=
\bigl(F_{f(\lambda_t,N_t)}^{-1}(U_{t+1}),f(\lambda_t,N_t)\bigr),
\]
where
$F_\lambda$ denotes the cumulative distribution function of a
$\operatorname{Poisson}(\lambda)$ distribution. This gives us a representation
of $((N_t,\lambda_t))_{t\in\N}$ as a randomly perturbed dynamical
system with independent and identically distributed innovations. In
such a context and under a contractive condition similar to our
(\ref{2.3}), Diaconis and Freedman \cite{DF99} also proved existence
and uniqueness of a stationary distribution. To this end, these authors
used backward iterations to identify a random variable which has the
desired stationary distribution. The approach used here is more direct
and uses also elements of standard Markov chain theory as described~in
Meyn and Tweedie \cite{MT93}. Finally, we would like to mention that
Lasota and Mackey \cite{LM89} also proved the existence of a unique
stationary distribution under conditions similar to our (\ref{2.3}) and
(\ref{2.4}); see, in particular, equations~(2) and~(3) in their paper.
Their proof contains similar ingredients to our proof; however, it is
more analytic in nature while we establish a coupling to represent the
convergence facts in a simple stochastic language.
\end{rem}

\section{Absolute regularity of the count process and ergodicity}
\label{S3}

In this section, we state the main results of our paper, absolute
regularity of the count process and, as a consequence, ergodicity of
the bivariate process $((N_t,\lambda_t))_t$. Actually, Grunwald,
Hyndman, Tedesco and Tweedie \cite{GHTT00}, Case~II of Proposition~3,
Streett \cite{Str00} and Davis, Dunsmuir and Streett \cite{DDS03}
proved ergodicity in special cases. However, they made heavy use of the
particular form of their link function~$f$ and could show that
Doeblin's condition is fulfilled. Hence, they could employ Markov chain
technology to prove ergodicity. We cannot use this approach in the case
considered here since Doeblin's condition will not be satisfied in
general. Another commonly used approach to proving ergodicity, which is
not restricted to the case of Markov chains, consists in proving first
strong mixing as a sufficient condition for ergodicity. It turns out,
however, that the bivariate process $((N_t,\lambda_t))_t$ is not
strongly mixing in general; a counterexample is given in
Remark~\ref{R3} below. The problem lies in the discreteness of the
distribution of the ``innovations'' $N_t$ while the $\lambda_t$ take
values on a continuous scale. This makes the commonly used coupling
approach to proving mixing properties of Markov chains impossible. To
give some idea why a discrete distribution of the innovations may cause
problems, we recall the well-known example of a stationary AR(1)
process, $X_t=\theta X_{t-1}+\varepsilon_t$, where the innovations
are independent with $P(\varepsilon_t=1)=P(\varepsilon_t=-1)=1/2$
and $0<|\theta|\leq1/2$. This process has a stationary
distribution supported on $[-2,2]$. It follows from the above model
equation that $X_t$ has, with probability~1, the same sign as
$\varepsilon_t$. Hence, we could perfectly recover
$X_{t-1},X_{t-2},\ldots$ from $X_t$, which clearly excludes any of the
common mixing properties. (Rosenblatt \cite{Ros80} mentioned the fact
that a process similar to $(X_t)_{t\in\Z}$ is purely deterministic
going backwards in time. A rigorous proof that it is not strong mixing
was given by Andrews \cite{And84}.) On the other hand, we can prove
absolute regularity for the (univariate) count process $(N_t)_t$. For
this purpose, the discrete nature of the distribution of the~$N_t$ does
not harm. To see why, note that we have either $\pi(\{0,0\})=1$ or
$P(\lambda_{t-1}>0 \mbox{ or } \lambda_t>0)=1$; see the proof of
part~(iii) of Theorem~\ref{T2.1}. Therefore, the support of the
conditional distribution of $N_{t+2}$ given $N_t,N_{t-1},\ldots$ is
equal to the support of the stationary distribution of the $N_t$ and we
can actually construct a successful coupling. Since absolute regularity
implies strong mixing, we immediately obtain ergodicity of the count
process $(N_t)_t$. Moreover, as a by-product of our coupling, we see
that the random intensities~$\lambda_t$ can be expressed as measurable
functionals of past variables of the count process. Hence, we finally
obtain the desired ergodicity of the bivariate process
$((N_t,\lambda_t))_t$.

It was stated in Section~\ref{S2} that the bivariate process
$((N_t,\lambda_t))_t$ has a unique stationary distribution under the
contractive condition (\ref{2.3}). In this section, we will assume
throughout that this process is in its stationary regime. Moreover, it
proves to be quite convenient to have a two-sided stationary version,
with time domain~$\Z$ rather than~$\N$, which exists by Kolmogorov's
extension theorem; see~Durrett \cite{Dur91}, page~293. Here is the main
result of the paper.

\begin{tthhmm}\label{T3.1}
Suppose that the bivariate chain
$((N_t,\lambda_t))_{t\in\Z}$ is in its stationary regime and obeys~(\ref{2.1})--(\ref{2.3}). Then
\begin{longlist}[(iii)]
\item The count process $(N_t)_{t\in\Z}$ is absolutely regular
with coefficients satisfying
\[
\beta(n) \leq2 E\lambda_1
\kappa^{n-1}/(1-\kappa_1).
\]
\item There exists a measurable
function $g\dvtx \N_0^{\infty}:=\{(n_1,n_2,\ldots)\dvt
n_i\in\N_0\}\longrightarrow[0,\infty)$ such that
$\lambda_t=g(N_{t-1},N_{t-2},\ldots)$ holds almost surely.
\item The process $((N_t,\lambda_t))_{t\in\Z}$ is ergodic.
\item $E \lambda_1^2 <\infty$.
\end{longlist}
\end{tthhmm}

\begin{rem}
In the case of a so-called $\operatorname{INGARCH}(1,1)$ process where $\lambda_t$
is specified as $\lambda_t=\theta_0+\theta_1\lambda_{t-1}+\theta_2
N_{t-1}$, Ferland, Latour and Oraichi \cite{FLO06} proved the
stronger result that all moments of $\lambda_t$ and $N_t$ are finite.
Since it follows from (\ref{2.4}) that $\lambda_t\leq
f(0,0)+\kappa_1\lambda_{t-1}+\kappa_2 N_{t-1}$, we conjecture that
their result can be generalized by simple majorization arguments to our
more general framework. However, since higher-than-second moments are
not needed for the purposes of this paper, we do not make the attempt
to adapt their proof, which was already quite involved in the special
case
of a linear specification of $\lambda_t$.
\end{rem}

\begin{rem}\label{R3}
Theorem~\ref{T3.1} states that the count process
$(N_t)_{t\in\Z}$ is absolutely regular and, therefore, also strongly
mixing under condition (\ref{2.3}). This allows us to conclude that the
bivariate process $((N_t,\lambda_t))_{t\in\Z}$ is ergodic. However, the
process $((N_t,\lambda_t))_{t\in\Z}$ and even the intensity process
$(\lambda_t)_{t\in\Z}$ alone are not strongly mixing in general. To see
this, consider the specification $f(\lambda,y)=g(\lambda)+y/2$,
where~$g$ is strictly monotone and satisfies $0<c_1\leq
g(\lambda)<0.5$ and $|g(\lambda)-g(\lambda')|\leq
c_2|\lambda-\lambda'|$ for some $c_2<0.5$ and for all
$\lambda,\lambda'$. Then~$f$ satisfies our contractive condition
(\ref{2.3}). Using the fact that $g(\lambda)\in[c_1,0.5),$ we
obtain that $2g(\lambda_{t-1})=2\lambda_t-[2\lambda_t],$ which
implies that we can perfectly recover $\lambda_{t-1}$ once we know the
value of $\lambda_t$. Iterating this argument, we see that we can
recover from $\lambda_t$ the entire past of the intensity process.
Taking into account that the above choice of~$f$ excludes the case that
the intensity process is purely non-random, we conclude that a
stationary version of $(\lambda_t)_{t\in\Z}$
cannot be strongly mixing.
\end{rem}

\begin{rem}
The primary intention of the author was to devise a method of
proving ergodicity of certain count processes. This is done, mainly for
clarity of presentation, for the simple case where the intensity
depends only on one lagged value of the count process and the intensity
process. In contrast to previous work in this area, the coupling
approach used here does not require Markovianity of the process. The
results of this paper, and in particular the ergodicity stated in
Theorem~\ref{T3.1}, can be generalized to more complex models with more
than one or even infinitely many lagged variables. Moreover, it seems
to be possible to include covariates, at least if they are exogeneous.
These generalizations are well beyond the scope
of this paper and should be the subject of future research.
\end{rem}

\section{A specification test for the intensity function}\label{S4}

There might be good reasons for assuming that the count variables are
Poisson distributed, conditioned on the past. However, a particular
specification for the intensity function seems to be more questionable
and such a choice should be supported by a statistical test. Here we
propose a test statistic that was originally designed for testing
overdispersion in the context of i.i.d.~observations; see Lee
\cite{Lee86} and Cameron and Trivedi \cite{CT86}.

Assume that we have observations $N_1,\ldots,N_n$ from a stationary
process $((N_t,\lambda_t))_{t\in\Z}$ obeying~(\ref{2.1}) and
(\ref{2.2}) and that we want to test the simple hypothesis
\[
H_0\dvt
f= f_0 \quad\mbox{against}\quad H_1\dvt f \neq
f_0,
\]
for some $f_0$ satisfying (\ref{2.3}), or the composite
hypothesis
\[
H_0'\dvt f\in\{f_\theta\dvt \theta\in\Theta\}
\quad\mbox{against}\quad H_1'\dvt f\notin\{f_\theta\dvt
\theta\in\Theta\},
\]
where $\Theta\subseteq\R^d$ and the
$f_\theta$
satisfy (\ref{2.3}).

To motivate a particular test statistic, pretend that we additionally
observe the starting value $\lambda_1$ of the intensity process. Then
we could take, for testing $H_0$ against $H_1$, the statistic
\[
T_{n,0} = \frac{1}{\sqrt{n}} \sum_{t=1}^n \{ (N_t -
\lambda_t^0)^2 - N_t \},
\]
where $\lambda_1^0=\lambda_1$
and, for $t=2,\ldots,n$, the $\lambda_t^0$ are recursively defined as
$\lambda_t^0=f_0(\lambda_{t-1}^0,N_{t-1})$. The idea behind this
statistic is very simple. If the the intensity function~$f$ is
correctly specified, then $\lambda_t^0=\lambda_t$, which implies
$E[ (N_t-\lambda_t^0)^2 - N_t ] = 0$ and, as stated
in Proposition~\ref{P4.2} below,
$T_{n,0}\stackrel{d}{\longrightarrow}{\mathcal N}(0,
2E\lambda_1^2)$. On the other hand, if~$f$ is not correctly specified
by~$f_0$, then the random variables $(N_t-\lambda_t)^2-N_t$ are not
centered and we can expect consistency of the test.

In the more relevant case of unknown $\lambda_1$, we replace this by
any arbitrarily chosen, random or non-random, starting value
$\widetilde{\lambda}_1$, then define recursively
$\widetilde{\lambda}_t=f_0(\widetilde{\lambda}_{t-1},N_{t-1})$, for
$t=2,\ldots,n$, and take the test statistic
\[
T_n =
\frac{1}{\sqrt{n}} \sum_{t=1}^n \{ (N_t - \widetilde{\lambda}_t)^2
- N_t \}.
\]
In the case of testing $H_0'$ against $H_1'$,
we estimate the parameter~$\theta$ by some estimator
$\widehat{\theta}_n$ first and take then the test statistic
\[
\widehat{T}_n = \frac{1}{\sqrt{n}} \sum_{t=1}^n \{ (N_t -
\widehat{\lambda}_t)^2 - N_t \}.
\]
Here
$\widehat{\lambda}_1$ is again any starting value and
$\widehat{\lambda}_t=f_{\widehat{\theta}_n}(\widehat{\lambda
}_{t-1},N_{t-1})$,
for $t=2,\ldots,n$.

\begin{rem}
In the context of independent observations, Lee \cite{Lee86} and
Cameron and Trivedi \cite{CT86} considered a test statistic similar to
ours for testing the Poisson hypothesis against the alternative that
the distribution belongs to the so-called Katz family of distributions.
This family contains as special cases the Poisson, negative binomial
and binomial distributions. While the variance equals the mean in the
Poisson case, the latter two classes contain distributions for which
the variance mean ratio is strictly greater or less than~1,
respectively. Therefore, Lee \cite{Lee86} and Cameron and Trivedi
\cite{CT86} interpreted their tests as tests for over- or
underdispersion. The same test statistic was also suggested in Cox
\cite{Cox83}. It was also used by Br\"{a}nn\"{a}s and Johansson
\cite{BJ94} for testing for the existence of a latent process in the
context of Poisson count models. Again, in the case of independent
data, Dean and Lawless \cite{DL89} and Dean \cite{Dea92} came up with
adjusted versions of Lee's and Cameron and Trivedi's test statistic
that have the same limit distribution as the unadjusted statistic but
are closer to this limit in small samples.
\end{rem}

We will prove that the above statistics, $T_{n,0}$, $T_n$ and
$\widehat{T}_n$, are asymptotically normal with the same limit. This
can be most easily done for $T_{n,0}$ since this statistic is a sum of
martingale differences that allows us to apply an appropriate central
limit theorem.

\begin{prop}\label{P4.2}
Suppose that the bivariate process is stationary
and obeys (\ref{2.1}) and (\ref{2.2}). If $H_0$ is true and $f_0$
satisfies the contractive condition (\ref{2.3}), then
\[
T_{n,0} \stackrel{d}{\longrightarrow} {\mathcal N}(0, 2E\lambda
_1^2).
\]
\end{prop}

Next, we will show that $T_n$ and $\widehat{T}_n$ have the same limit
distribution as $T_{n,0}$. To this end, we will simply show that the
difference between the former statistics to $T_{n,0}$ converges to zero
in probability. This is not surprising at all for $T_n$ since it
follows from (\ref{2.3}) that
$|\widetilde{\lambda}_t-\lambda_t|\leq\kappa_1^{t-1}|\widetilde
{\lambda}_1-\lambda_1|$.
The following lemma shows that $\widehat{\lambda}_t$ will also be
close to $\lambda_t$ if $\widehat{\theta}_n$ is a $\sqrt{n}$-consistent
estimator of~$\theta$ and if $f_\theta(\lambda,y)$ is a smooth function
in~$\theta$.

\begin{lem}\label{L4.3}
Suppose that the bivariate process is stationary and
obeys (\ref{2.1}) and (\ref{2.2}) with $f=f_{\theta_0}$. We assume that
$\widehat{\theta}_n-\theta_0=\mathrm{O}_P(n^{-1/2})$. Furthermore, we assume
that there exist $C<\infty$, $\kappa_1,\kappa_2\geq0$ with
$\kappa:=\kappa_1+\kappa_2<1$ such that
\begin{longlist}[(ii)]
\item $ |f_{\theta'}(\lambda,y) -
f_{\theta_0}(\lambda,y)| \leq C \|\theta' - \theta_0\|
(\lambda+ y + 1)$ $\forall\lambda, y,$
\item $ |f_{\theta'}(\lambda,y) -
f_{\theta'}(\widetilde{\lambda},\widetilde{y})| \leq\kappa_1
|\lambda- \widetilde{\lambda}| + \kappa_2 |y - \widetilde{y}|$
\end{longlist}
hold for all $\theta'\in\Theta$ with $\|\theta'-\theta_0\|\leq
\delta$,
for some $\delta>0$.

Then
\[
\sum_{t=1}^n (\lambda_t - \widehat{\lambda}_t)^2 = \mathrm{O}_P(1).
\]
\end{lem}

We think that the above assumption on the estimator
$\widehat{\theta}_n$ is a realistic one in many cases. It is fulfilled,
for example, by the conditional maximum likelihood estimator studied in
Fokianos, Rahbek and Tj{\o}stheim \cite{FRT09}.

\begin{tthhmm}\label{T4.4}
Suppose that the assumptions of Lemma~\ref{L4.3} are
fulfilled.

Then
\[
\widehat{T}_n \stackrel{d}{\longrightarrow} {\mathcal N}(0, 2
E\lambda_1^2).
\]
\end{tthhmm}

\begin{rem}
The same assertion holds true for $T_n$ instead of
$\widehat{T}_n$ since this
is obviously a special case of that considered in Theorem~\ref{T4.4}.
\end{rem}

Note that the limit distribution of $\widehat{T}_n$ still contains the
parameter $E\lambda_1^2$ that is usually not known in advance and has
to be estimated. We obtain from Lemma~\ref{L4.3} by the Minkowski
inequality that $| \sqrt{ n^{-1} \sum_{t=1}^n
\widehat{\lambda}_t^2 } - \sqrt{ n^{-1} \sum_{t=1}^n \lambda_t^2 }
| \leq\sqrt{ n^{-1} \sum_{t=1}^n (\lambda_t -
\widehat{\lambda}_t)^2 } = \mathrm{O}_P( n^{-1/2} )$, which leads in
conjunction with ergodicity of $(\lambda_t)_{t\in\Z}$ to
%
\begin{equation}
\label{4.1} \frac{1}{n} \sum_{t=1}^n \widehat{\lambda}_t^2
\stackrel{P}{\longrightarrow} E \lambda_1^2.
\end{equation}
For a prescribed size~$\alpha\in(0,1)$, we propose a test for $H_0'$
against $H_1'$ as
\[
\varphi_n = I\Biggl( \Biggl( (2/n) \sum_{t=1}^n
\widehat{\lambda}_t^2 \Biggr)^{-1/2} \widehat{T}_n > u_\alpha
\Biggr),
\]
where $u_\alpha=\Phi^{-1}(1-\alpha)$ denotes the
$(1-\alpha)$-quantile of the standard normal distribution. From
Theorem~\ref{T4.4} and (\ref{4.1}) we conclude that this test has
asymptotically the correct size.

\begin{tthhmm}\label{T4.5}
Suppose that the assumptions of Lemma~\ref{L4.3} are
fulfilled and that $f_{\theta_0}(0,0)>0$. Then we have, under~$H_0'$,
\[
\Biggl( (2/n) \sum_{t=1}^n \widehat{\lambda}_t^2 \Biggr)^{-1/2}
\widehat{T}_n \stackrel{d}{\longrightarrow} {\mathcal N}(0,1),
\]
which implies that
\[
P( \varphi_n = 1 ) \mathop{\longrightarrow}_{n\to\infty} \alpha.
\]
\end{tthhmm}

\section{Proofs}\label{SP}

As already mentioned in the text, the main results of this paper,
Theorems~\ref{T2.1} and \ref{T3.1}, are both proved by coupling
arguments. Necessary technical prerequisites are summarized in the
following lemma.

\begin{lem}\label{LA.1}
For arbitrary $\lambda_1,\lambda_2\geq0$, we can
construct on an appropriate probability space $X_1\sim
\operatorname{Poisson}(\lambda_1)$ and $X_2\sim\operatorname{Poisson}(\lambda_2)$ such
that
\begin{longlist}[(ii)]
\item $ E|X_1-X_2| = |\lambda_1-\lambda_2|$,
\item $ P(X_1\neq X_2) \leq|\lambda_1-\lambda_2|$.
\end{longlist}
\end{lem}

\begin{pf}
Let, without loss of generality, $\lambda_1\leq\lambda_2$. We take
independent random variables $X_1\sim\operatorname{Poisson}(\lambda_1)$,
$Z\sim\operatorname{Poisson}(\lambda_2-\lambda_1)$ and define $X_2=X_1+Z$.
Then $X_2\sim\operatorname{Poisson}(\lambda_2)$,
\[
E|X_1-X_2| = EZ =|\lambda_1-\lambda_2|
\]
and
\[
P(X_1\neq X_2) = P(Z\neq0)\leq EZ = |\lambda_1-\lambda_2|.
\]
\end{pf}

\begin{pf*}{Proof of Theorem~\ref{T2.1}}
As mentioned above, we could use the fact that
$((N_t,\lambda_t))_{t\in\N}$ is a weak Feller chain that is bounded in
probability on average to conclude from Theorem~12.1.2(ii) in Meyn and
Tweedie \cite{MT93} that it has at least one stationary distribution.
Uniqueness could then eventually be derived from the contraction
property (\ref{2.3}). We think, however, that it is more instructive
for the reader when a self-contained proof that uses arguments closely
tied to the particular case at hand is presented.

Let $P_\lambda^t$ be the conditional distribution of $(N_t,\lambda_t)$
given $\lambda_1=\lambda$, where $\lambda\in[0,\infty)$ is an
arbitrarily chosen but fixed starting value. It follows from
(\ref{2.5}) that the sequence of distributions $(P_\lambda^t)_{t\in
\N}$
is tight. Hence, there exists a subsequence $(n_k)_{k\in\N}$ of $\N$
such that $P_\lambda^{n_k}$ converges weakly to some probability
measure $\pi_\lambda$, as $k\to\infty$. We will show that this limit
does not depend on the starting value~$\lambda$ and that the full
sequence $(P_\lambda^n)_{n\in\N}$ converges. This will immediately
imply that $\pi_\lambda$ is a stationary distribution that is unique.

The latter conclusions will follow after we have derived a few
convergence properties of the process. To this end, we construct, on an
appropriate probability space $(\Omega',{\mathcal A}',P')$, two Markov
chains $((N_t',\lambda_t'))_{t\in\N}$ and
$((N_t'',\lambda_t''))_{t\in\N}$ with transition laws according
to (\ref{2.1}) and (\ref{2.2}) and with starting values $\lambda_1'$
and $\lambda_1''$, respectively. We construct these chains
iteratively. Given $\lambda_1'$ and $\lambda_1''$, (i) of
Lemma~\ref{LA.1} allows us to construct $N_1'$ and $N_1''$ in such a
way that
\[
E(| N_1' - N_1'' | | \lambda_1',
\lambda_1'' ) = | \lambda_1' - \lambda_1'' |.
\]
The
values of $\lambda_2'$ and $\lambda_2''$ are then given by equation
(\ref{2.2}) and it follows from (\ref{2.3}) that
\begin{eqnarray*}
E(|\lambda_2' - \lambda_2''| | \lambda_1', \lambda_1'')
&\leq&
\kappa_1 |\lambda_1' - \lambda_1''| + \kappa_2E(|N_1' - N_1''| | \lambda_1', \lambda_1'' )
\\
&=& \kappa| \lambda_1' - \lambda_1'' |.
\end{eqnarray*}
In the next step we
can construct $N_2'$ and $N_2''$ such that $E(|N_2'-N_2''| |
\lambda_1', \lambda_1'', \lambda_2', \lambda_2'')
=|\lambda_2'-\lambda_2''|,$ which also implies that
\[
E(|N_2' - N_2''| | \lambda_1', \lambda_1'' )\leq\kappa| \lambda_1' - \lambda_1'' |.
\]
Now we can proceed
in the same way and construct the pairs $(N_3',N_3''),
(N_4',N_4''),\ldots.$ With the above construction, we obtain that
%
\begin{equation}
\label{pt21.1}
E(| \lambda_t' - \lambda_t'' | |\lambda_1', \lambda_1'' ) \leq\kappa^{t-1} | \lambda_1' -\lambda_1'' |.
\end{equation}
and
%
\begin{equation}
\label{pt21.2}
E(| N_t' - N_t'' | | \lambda_1',\lambda_1'' ) \leq\kappa^{t-1} | \lambda_1' -\lambda_1'' |.
\end{equation}
Hence, it follows that $(P_{\lambda_1'}^{n_k})_{k\in\N}$ and
$(P_{\lambda_1''}^{n_k})_{k\in\N}$ converge for any choice of
$\lambda_1'$ and $\lambda_1''$ to the same limit, which we denote by
$\pi$ in the following. Now we can translate this result to a
convergence result for the conditional distributions of the Markov
chain $((N_t,\lambda_t))_{t\in\N}$. Since the above convergence is
uniform in~$\lambda_1'$ over compact sets and since~$f$ as a continuous
function maps compact subsets of $[0,\infty)\times\N_0$ to compact
subsets of $[0,\infty),$ we obtain that
%
\begin{equation}\label{pt21.3}
\sup_{x\in K} | P^{n_k}(x,A) - \pi(A) |\mathop{\longrightarrow}_{k\to\infty} 0
\end{equation}
holds for every compact subset~$K$ of $\N_0\times[0,\infty)$ and every
$\pi$-continuity set~$A$. Here $P^{n}(x,A)=P((N_{t+n},\lambda_{t+n})\in
A|(N_t,\lambda_t)=x)$ denotes the $n$-step transition probability
of the bivariate process. Equation (\ref{pt21.3}) will allow us to show
convergence of the full sequence. For any $n\in\N$, let $k(n)$ be the
largest integer such that $n_{k(n)}< n$. From tightness of
$(P_\lambda^n)_{n\in\N}$ and (\ref{pt21.3}), we conclude that
%
\begin{equation}
\label{pt21.4} P_\lambda^n = \int P^{n_{k(n)}}(x, \cdot)
P_\lambda^{n-n_{k(n)}}(\mathrm{d}x) \quad\Longrightarrow\quad\pi\qquad\mbox{for all
} \lambda\in[0,\infty).
\end{equation}
It follows directly from this equation that~$\pi$ is a stationary
distribution. To see this, observe that it follows from (\ref{pt21.4})
that $Q_x^n:=n^{-1} \sum_{t=1}^n P_\lambda^t$ converges weakly
to~$\pi$. Furthermore, it also follows that
$\widetilde{Q}_\lambda^n:=n^{-1}\sum_{t=1}^n
P_\lambda^{t+1}\Longrightarrow\pi$, that is,
%
\begin{equation}
\label{pt21.5}
\widetilde{Q}_\lambda^n(A)
\mathop{\longrightarrow}_{k\to\infty} \pi(A)
\end{equation}
if $A$ is a $\pi$-continuity set, that is, $\pi(\partial A)=0$. If~$A$
is an open set, then $x\mapsto P^1(x,A)$ is a continuous and bounded
function. Therefore,
%
\begin{equation}
\label{pt21.6} \widetilde{Q}_\lambda^n(A) = \int P^1(x,A)
Q_\lambda^{n}(\mathrm{d}x) \mathop{\longrightarrow}_{n\to\infty} \int
P^1(x,A) \pi(\mathrm{d}x).
\end{equation}
From (\ref{pt21.5}) and (\ref{pt21.6}) we obtain that the probability
measures~$\pi$ and $\int P^1(x,\cdot) \pi(\mathrm{d}x)$ coincide for all open
$\pi$-continuity sets~$A$. Since these sets are stable under finite
intersections and generate $2^{\N_0}\otimes{\mathcal B,}$ we conclude
that
\[
\pi(A) = \int P^1(x,A) \pi(\mathrm{d}x)\qquad \forall A\in
2^{\N_0}\otimes{\mathcal B},
\]
that is, $\pi$ is actually a
stationary distribution. Let $\pi'$ be an arbitrary distribution. Then
we obtain by majorized convergence, for any $\pi$-continuity set~$A$,
\[
\int P^n(x,A) \pi'(\mathrm{d}x)
\mathop{\longrightarrow}_{n\to\infty} \int\pi(A)
\pi'(\mathrm{d}x) = \pi(A).
\]
If $\pi'$ is a stationary distribution, then
we also have that $\int P^n(x,A) \pi'(\mathrm{d}x)=\pi'(A)$, which
implies that $\pi=\pi'$. Hence, (i) is proved.

We obtain from (\ref{2.5}) and by Theorem~5.3 in Billingsley
\cite{Bil68} that
\begin{displaymath}E_\pi\lambda_1 \leq\liminf_{t\to\infty}
E(\lambda_t|\lambda_1=0) \leq f(0,0) / (1-\kappa),
\end{displaymath}
which
proves (ii).

To see (iii), note that $f(0,0)=0$ implies by (\ref{2.5}) that
$E(\lambda_t|\lambda_1=0)=0$ holds for all~$t$ which in turn
implies that $\pi(\{0,0\})=1$. On the other hand, if $f(0,0)>0$, then
we can conclude that $P(\lambda_{t-1}>0 \mbox{ or } \lambda_t>0) =
P(\lambda_{t-1}>0) + P(\lambda_{t-1}=0, \lambda_t>0) =
P(\lambda_{t-1}>0) + P(\lambda_{t-1}=0) = 1$. This implies that
$((N_t,\lambda_t))_{t\in\Z}$ cannot be non-random, as required.
\end{pf*}

\begin{pf*}{Proof of Theorem~\ref{T3.1}}
Let, for $-\infty\leq k\leq l\leq\infty$, ${\mathcal
B}^N_{k,l}=\sigma(N_k,\ldots,N_l)$. Recall that the coefficients of
absolute regularity of the count process $(N_t)_{t\in\N}$ are defined
as
\[
\beta(n) = E\Bigl[ \sup_{A\in{\mathcal B}^N_{n,\infty
}} \bigl|
P(A|{\mathcal B}^N_{-\infty,0}) - P(A) \bigr| \Bigr].
\]
Hence,
\[
\beta(n) \leq E\Bigl[ \sup_{A\in{\mathcal
B}^N_{n,\infty}} \bigl|
P(A|\sigma(\lambda_1,N_0,N_{-1},\ldots)) - P(A) \bigr| \Bigr].
\]
Furthermore, it follows from $(N_{n},N_{n+1},\ldots)|
\sigma(\lambda_1,N_0,N_{-1},\ldots) = (N_{n},N_{n+1},\ldots)|
\sigma(\lambda_1)$ that
%
\begin{equation}
\label{pt31.1}
\beta(n) \leq E\Bigl[ \sup_{A\in{\mathcal
B}^N_{n,\infty}} \bigl| P(A|\sigma(\lambda_1)) - P(A) \bigr| \Bigr].
\end{equation}
Let ${\mathcal B}^\infty$ be the $\sigma$-field in
$\R^\infty=\{(x_1,x_2,\ldots)\dvt x_i\in\R\}$ generated by the cylinder
sets, that is,
\[
{\mathcal B}^\infty= \sigma( \{
B\times\R^\infty\dvt B\in{\mathcal B}^k, k\in\N\} ).
\]
We can
rewrite (\ref{pt31.1}) in terms of the process variables as
%
\begin{equation}
\label{pt31.2}
\beta(n) \leq E\biggl[ \sup_{A\in{\mathcal
B}^\infty} \bigl| P\bigl((N_n,N_{n+1},\ldots)\in A|\lambda_1 \bigr) -
P\bigl((N_n,N_{n+1},\ldots)\in A \bigr) \bigr| \biggr].
\end{equation}
We will derive an upper estimate for the right-hand side of
(\ref{pt31.2}) via a coupling approach similar to that in the proof of
Theorem~\ref{T2.1}. To this end, we will construct on an appropriate
probability space $(\Omega',{\mathcal A}',P')$ two versions of the
bivariate process, $((N_t',\lambda_t'))_{t\in\N}$ and
$((N_t'',\lambda_t''))_{t\in\N}$, where the starting values
$\lambda_1'$ and $\lambda_1''$ are independent and distributed
according to the stationary law~$\pi$. Since, for any $A\in{\mathcal
B}^\infty$,
\[
P\bigl( (N_n'',N_{n+1}'',\ldots)\in A|
\lambda_1' \bigr) = P\bigl( (N_n,N_{n+1},\ldots)\in A \bigr)
\]
it follows that
\begin{eqnarray*}
&&
\bigl| P\bigl( (N_n,N_{n+1},\ldots)\in A |\lambda_1=u \bigr) -P\bigl( (N_n,N_{n+1},\ldots)\in A \bigr) \bigr|
\\
&&\quad=
\bigl| P\bigl( (N_n',N_{n+1}',\ldots)\in A |\lambda_1'=u\bigr) -P\bigl( (N_n'',N_{n+1}'',\ldots)\in A |\lambda_1'=u \bigr) \bigr| \\
&&\quad\leq
P\bigl( (N_n',N_{n+1}',\ldots) \neq(N_n'',N_{n+1}'',\ldots) |\lambda_1'=u \bigr).
\end{eqnarray*}
Therefore, we obtain that
%
\begin{equation}\label{pt31.3}
\beta(n) \leq P\bigl(
(N_n',N_{n+1}',\ldots) \neq(N_n'',N_{n+1}'',\ldots) \bigr).
\end{equation}
Hence, to estimate $\beta(n)$, we will construct a coupling such that
the processes $(N_t')_{t\in\N}$ and $(N_t'')_{t\in\N}$ coalesce
after~$n$ steps with a high probability.

Using exactly the same construction as in the proof of
Theorem~\ref{T2.1}, we can successively construct pairs
$(N_1',N_1''),(N_2',N_2''),\ldots$ such that
\[
E(|\lambda_n' - \lambda_n''| | \lambda_1', \lambda_1''
) \leq\kappa^{n-1} |\lambda_1' - \lambda_1''|.
\]
From here on we deviate from the approach in the proof of
Theorem~\ref{T2.1}, where we constructed all pairs $(N_t',N_t'')$
such that their mean distance was small. By (ii) of Lemma~\ref{LA.1},
we can construct $N_n'$ and $N_n''$ such that
\[
P(N_n'\neq N_n'' | \lambda_1', \lambda_1'' ) \leq
\kappa^{n-1} |\lambda_1' - \lambda_1''|.
\]
If the event
$\{N_n' = N_n''\}$ occurs, then (\ref{2.3}) reduces to
\[
|
\lambda_{n+1}' - \lambda_{n+1}'' | \leq\kappa_1 |
\lambda_n' - \lambda_n'' |,
\]
which allows us to construct the
next pair $(N_{n+1}',N_{n+1}'')$ such that
\[
P( N_n'=
N_n'', N_{n+1}'\neq N_{n+1}'' | \lambda_1', \lambda_1''
) \leq\kappa_1 \kappa^{n-1} |\lambda_1' -
\lambda_1''|.
\]
Continuing in the same way, we arrive at
\[
P(N_n'= N_n'', \ldots, N_{n+k-1}'= N_{n+k-1}'',
N_{n+k}'\neq N_{n+k}'' | \lambda_1', \lambda_1'' )
\leq\kappa_1^k \kappa^{n-1} |\lambda_1' -
\lambda_1''|.
\]
Hence, we finally obtain that
%
\begin{eqnarray}
&&
P\bigl( (N_n',N_{n+1}',\ldots) \neq(N_n'',N_{n+1}'',\ldots)
\bigr)\nonumber
\\
&&\quad=
P( N_n' \neq N_n'' ) + \sum_{k=1}^\infty P( N_n'=N_n'',\ldots, N_{n+k-1}'= N_{n+k-1}'',N_{n+k}'\neq N_{n+k}'' )
\qquad \\
&&\quad\leq
C_0 \kappa^{n-1}/(1-\kappa_1),\nonumber
\end{eqnarray}
where
$C_0:=E|\lambda_1'-\lambda_1''|\leq2E\lambda_1<\infty$. This
yields, in conjunction with (\ref{pt31.3}), Assertion~(i).

To show~(ii), define the functions $f_1=f$ and, for $d\geq2$,
$f_d(\lambda;n_1,\ldots,n_d)=f_{d-1}(f(\lambda,n_d);\break n_1,\ldots,n_{d-1})$,
where $n_1,\ldots,n_d\in\N_0$ and $\lambda\geq0$. It is clear from
(\ref{2.2}) that
\[
\lambda_t =
f_d(\lambda_{t-d};N_{t-1},\ldots,N_{t-d}).
\]
It follows from
(\ref{2.3}) that
\[
E| \lambda_t -
f_d(0;N_{t-1},\ldots,N_{t-d}) | \leq\kappa_1^d E
\lambda_{t-d}.
\]
Hence, as $d\to\infty$,
$f_d(0;N_{t-1},\ldots,N_{t-d})$ converges in $L_1$ to~$\lambda_t$. By
taking an appropriate subsequence, we also get almost sure convergence.
This means that there exists a measurable function $f_\infty\dvtx
\N_0^{\infty}\longrightarrow[0,\infty)$ such that
%
\begin{equation}
\label{aslimit} \lambda_t = f_\infty(N_{t-1},N_{t-2},\ldots)\qquad
\mbox{almost surely}.
\end{equation}
By stationarity, (\ref{aslimit}) holds for all $t\in\Z$, which
proves~(ii).

To show (iii), we first recall the well-known fact that absolute
regularity implies strong mixing. That is, it follows from~(i) that
%
\begin{equation}
\label{pt31.4}
\alpha(n) = \sup_{A\in{\mathcal B}^N_{-\infty,0},
B\in{\mathcal B}^N_{n,\infty}} | P(A\cap B) - P(A) P(B) |
\mathop{\longrightarrow}_{n\to\infty} 0;
\end{equation}
see Doukhan \cite{Dou94}, page~20. Furthermore, strong mixing implies
ergodicity; see~Remark~2.6 on page~50 in combination with
Proposition~2.8 on page~51 in Bradley \cite{Bra07}. Finally, we
conclude from the representation (\ref{aslimit}) by
Proposition~2.10(ii) in Bradley \cite{Bra07}, page~54, that the bivariate
process $((N_t,\lambda_t))_{t\in\Z}$ is also ergodic.

To prove (iv), we study the asymptotics of the process
$((\widetilde{N}_t,\widetilde{\lambda}_t))_{t\in\N}$ obeying
(\ref{2.1}), (\ref{2.2}) and (\ref{2.3}), which is started with
$\widetilde{\lambda}_1\equiv0$. We obtain from (\ref{2.4}) and
$E(N_t^2|\lambda_t)=\lambda_t^2+\lambda_t$ that
\begin{eqnarray*}
E(\widetilde{\lambda}_t^2 |\widetilde{\lambda}_{t-1})
& \leq&
E\bigl( \bigl(f(0,0) + \kappa_1 \widetilde{\lambda}_{t-1} + \kappa_2\widetilde{N}_{t-1}\bigr)^2|\widetilde{\lambda}_{t-1} \bigr)
\\
& = &
\bigl( f(0,0) + \kappa\widetilde{\lambda}_{t-1} \bigr)^2+ \kappa_2^2 \widetilde{\lambda}_{t-1}
\\
& \leq&
K_0 + \bar{\kappa} \widetilde{\lambda}_{t-1}^2,
\end{eqnarray*}
for
any $\bar{\kappa}>\kappa$ and appropriate $K_0=K_0(\bar{\kappa})$. We
choose $\bar{\kappa}\in(\kappa,1)$. Then we obtain that
\[
E(\widetilde{\lambda}_3^2 |\widetilde{\lambda}_1) \leq K_0 +\bar{\kappa} E( \widetilde{\lambda}_2^2 |\widetilde{\lambda}_1)
\leq K_0 + \bar{\kappa} ( K_0 + \bar{\kappa}\widetilde{\lambda}_1^2).
\]
Continuing in the same way we arrive at
the inequality
\[
E( \widetilde{\lambda}_t^2 |\widetilde{\lambda}_1) \leq K_0 (1 + \bar{\kappa} +\cdots+ \bar{\kappa}^{t-2}).
\]
Since
$\widetilde{\lambda}_t\stackrel{d}{\longrightarrow} \lambda_1$, we
conclude from Theorem~5.3 in Billingsley \cite{Bil68} that
\[
E
\lambda_1^2 \leq\liminf_{t\to\infty} E\widetilde{\lambda}_t^2
\leq K_0/(1- \bar{\kappa}),
\]
which proves (iii).
\end{pf*}

\begin{pf*}{Proof of Proposition~\ref{P4.2}}
We will use the central limit theorem (CLT) for martingale difference
arrays given on page~171 in Pollard \cite{Pol84}. We define the
filtration $({\mathcal B}_t)_{t\in\N}$ with ${\mathcal
B}_t=\sigma(\lambda_1,N_1,\ldots,N_t)$, for $t=0,1,\ldots,$ and we set
$Z_t=(N_t-\lambda_t)^2-N_t$. Since $N_t|{\mathcal
B}_{t-1}\sim\operatorname{Poisson}(\lambda_t),$ we obtain that
\[
E( Z_t
|{\mathcal B}_{t-1} ) = 0
\]
and
\[
E( Z_t^2 |{\mathcal
B}_{t-1} ) = 2\lambda_t^2.
\]
Hence, it follows from the
ergodicity stated in (iii) of Theorem~\ref{T3.1} that
\[
\frac{1}{n}
\sum_{t=1}^n E( Z_t^2 |{\mathcal B}_{t-1} )
\stackrel{\mathrm{{a.s.}}}{\longrightarrow} 2E \lambda_1^2.
\]
It remains to
verify the conditional Lindeberg condition,
\[
n^{-1} \sum_{t=1}^n E\bigl(
Z_t^2 I\bigl( \bigl|Z_t/\sqrt{n}\bigr| > \epsilon|{\mathcal B}_{t-1} \bigr) \bigr)
\stackrel{P}{\longrightarrow} 0\qquad \forall\epsilon>0.
\]
We
have $E[ n^{-1} \sum_{t=1}^n E( Z_t^2 I( |Z_t/\sqrt{n}| >
\epsilon|{\mathcal B}_{t-1} ) ) ] = E[ Z_1^2 I(
|Z_1/\sqrt{n}| > \epsilon) ]$, which tends to zero as $n\to\infty$
since $E\lambda_1^2<\infty$ implies that $EZ_1^2<\infty$. Hence, the
conditional Lindeberg condition is also satisfied and the assertion
follows from the CLT mentioned above.
\end{pf*}

\begin{pf*}{Proof of Lemma~\ref{L4.3}}
Assume for the time being that $\|\widehat{\theta}_n-\theta_0\|\leq
\delta$, which allows us to conveniently exploit the smoothness
assumptions on $f_\theta$. Then we obtain that
\begin{eqnarray*}
|\widehat{\lambda}_2 - \lambda_2|
& \leq&
|
f_{\widehat{\theta}_n}(\widehat{\lambda}_1,N_1) -f_{\widehat{\theta}_n}(\lambda_1,N_1) |
+ | f_{\widehat{\theta}_n}(\lambda_1,N_1) - f_{\theta_0}(\lambda_1,N_1) |
\\
& \leq&
\kappa_1 |\widehat{\lambda}_1 - \lambda_1| + C\|\widehat{\theta}_n-\theta_0\| (\lambda_1+N_1+1)
\end{eqnarray*}
and
\begin{eqnarray*}
|\widehat{\lambda}_3 - \lambda_3|
& \leq&
\kappa_1|\widehat{\lambda}_2 - \lambda_2|+ C \|\widehat{\theta}_n-\theta_0\| (\lambda_2+N_2+1)
\\
& \leq&
C \|\widehat{\theta}_n-\theta_0\| \{(\lambda_2+N_2+1) + \kappa_1 (\lambda_1+N_1+1) \} +\kappa_1^2 |\widehat{\lambda}_1 - \lambda_1|.
\end{eqnarray*}
Continuing
in the same way, we arrive at
\begin{eqnarray}\label{pl43.1}
&&|\widehat{\lambda}_t - \lambda_t| \nonumber
\\
&&\quad\leq
C \|\widehat{\theta}_n-\theta_0\| \{(\lambda_{t-1}+N_{t-1}+1) \nonumber
\\[-8pt]\\[-8pt]
&&\hphantom{\quad\leq
C \|\widehat{\theta}_n-\theta_0\| \{}{}+ \kappa_1(\lambda_{t-2}+N_{t-2}+1) +\cdots+\kappa_1^{t-2} (\lambda_1+N_1+1) \}\nonumber
\\
&&\qquad{}+ \kappa_1^{t-1} |\widehat{\lambda}_1 - \lambda_1|,\nonumber
\end{eqnarray}
which yields that
\begin{eqnarray*}
&&(\widehat{\lambda}_t - \lambda_t)^2
\\
&&\quad\leq
2 C^2 \|\widehat{\theta}_n-\theta_0\|^2 \{(\lambda_{t-1}+N_{t-1}+1) + \kappa_1(\lambda_{t-2}+N_{t-2}+1) +\cdots+\kappa_1^{t-2} (\lambda_1+N_1+1) \}^2
\\
&&\qquad{}+ 2 \kappa_1^{2t-2} (\widehat{\lambda}_1 -\lambda_1)^2
\end{eqnarray*}
holds for all $t\geq2$. Hence, we obtain under
$\|\widehat{\theta}_n-\theta_0\|\leq\delta$ that
\[
\sum_{t=1}^n
(\widehat{\lambda}_t - \lambda_t)^2 \leq\frac{2}{1-
\kappa_1^2} \Biggl\{ (\widehat{\lambda}_1 - \lambda_1)^2 + C^2
\|\widehat{\theta}_n-\theta_0\|^2 \Biggl( \sum_{t=1}^{n-1}
(\lambda_t + N_t + 1) \Biggr)^2 \Biggr\}.
\]
The right-hand
side is bounded in probability, which proves the assertion.
\end{pf*}

\begin{pf*}{Proof of Theorem~\ref{T4.4}}
We show that the difference between the test statistic $\widehat{T}_n$
and $T_{n,0}$ tends to zero in probability. This will yield the
assertion by Proposition~\ref{P4.2}. We have that
%
\begin{equation}\label{pt44.1}
\widehat{T}_n - T_{n,0} = \frac{1}{\sqrt{n}}\sum_{t=1}^n (\widehat{\lambda}_t - \lambda_t)^2 +
\frac{2}{\sqrt{n}} \sum_{t=1}^n (N_t-\lambda_t) (\lambda_t -\widehat{\lambda}_t).
\end{equation}
According to Lemma~\ref{L4.3}, the first term on the right-hand side
converges to zero in probability. The estimation of the second one,
however, is more delicate since $\widehat{\lambda}_t$ depends via
$\widehat{\theta}_n$ on the whole sample, which means that this term is
not a sum of martingale differences. To proceed, we take first any {\em
non-random} $\theta'$ with $\|\theta'-\theta_0\|\leq\delta$ and
consider the intensity process given by
$\lambda_1'=\widehat{\lambda}_1$ and, for $t=2,\ldots,n$,
$\lambda_t'=f_{\theta'}(\lambda_{t-1}',N_{t-1})$. We obtain in
complete analogy to (\ref{pl43.1}) that
\begin{eqnarray*}
&&|\lambda_t' - \lambda_t|
\\
&&\quad\leq
C \|\theta'-\theta_0\| \{ (\lambda_{t-1}+N_{t-1}+1)+ \kappa_1(\lambda_{t-2}+N_{t-2}+1)+ \cdots+ \kappa_1^{t-2} (\lambda_1+N_1+1) \}
\\
&&\qquad{}+
\kappa_1^{t-1} |\widehat{\lambda}_1 - \lambda_1|.
\end{eqnarray*}
Therefore, we obtain that
%
\begin{equation}
\label{pt44.2}
E\Biggl[ \Biggl| \frac{1}{\sqrt{n}} \sum_{t=1}^n
(N_t-\lambda_t)(\lambda_t-\lambda_t') \Biggr| I(|\widehat{\lambda}_1 -
\lambda_1 |\leq M) \Biggr] = \mathrm{O}( \| \theta' - \theta_0 \| +
n^{-1/2} ).
\end{equation}
Since $\|\widehat{\theta}_n-\theta_0\|=\mathrm{O}_P(n^{-1/2})$ it suffices
to establish (\ref{pt44.2}) on a sequence of grids ${\mathcal G}_n$ on
the set $\{\theta\in\Theta\dvt \|\theta-\theta_0\|\leq\epsilon_n^{-1}
n^{-1/2}\}$, where $\mbox{mesh}({\mathcal G}_n)\leq
\epsilon_nn^{-1/2}$, $\#{\mathcal G}_n\leq\epsilon_n n^{1/2}$,
for some null sequence $(\epsilon_n)_{n\in\N}$. It follows from
(\ref{pt44.2}) that
%
\begin{equation}
\label{pt44.3} \sup_{\theta'\in{\mathcal G}_n} \Biggl|
\frac{1}{\sqrt{n}} \sum_{t=1}^n (N_t-\lambda_t)(\lambda_t-\lambda_t')
\Biggr| = \mathrm{O}_P( \epsilon_n ).
\end{equation}
Moreover, for any value of $\widehat{\theta}_n$ with
$\|\widehat{\theta}_n-\theta_0\|\leq\epsilon_n^{-1}n^{-1/2}$ we
will find some $\theta'\in{\mathcal G}_n$ with
$\|\widehat{\theta}_n-\theta'\|\leq\epsilon_n n^{-1/2}$. Since
\begin{eqnarray*}
\Biggl| \frac{1}{\sqrt{n}} \sum_{t=1}^n
(N_t-\lambda_t)(\lambda'-\widehat{\lambda}_t) \Biggr| \leq\sqrt{
\frac{1}{n} \sum_{t=1}^n (N_t-\lambda_t)^2 } \sqrt{ \sum_{t=1}^n
(\lambda_t' - \widehat{\lambda}_t)^2 } = \mathrm{o}_P(1),
\end{eqnarray*}
we obtain,
in conjunction with (\ref{pt44.3}), that the second term on the
right-hand side of (\ref{pt44.1}) is $\mathrm{o}_P(1)$. This completes the
proof.
\end{pf*}

\section*{Acknowledgements}
This work was supported by the German Research Foundation DFG, project
NE 606/2-1. I thank Richard A.~Davis and an anonymous referee for
helpful comments on an earlier version of this paper. I also thank Paul
Doukhan and Konstantinos Fokianos for stimulating discussions.

%

\printhistory


\begin{thebibliography}{29}

\bibitem{And84}
%
\begin{barticle}[mr]
\bauthor{\bsnm{Andrews},~\bfnm{Donald W.~K.}\binits{D.W.K.}}
(\byear{1984}).
\btitle{Nonstrong mixing autoregressive processes}.
\bjournal{J. Appl. Probab.}
\bvolume{21}
\bpages{930--934}.
\bid{issn={0021-9002}, mr={0766830}}
\end{barticle}
%
\endbibitem

\bibitem{Bil68}
%
\begin{bbook}[mr]
\bauthor{\bsnm{Billingsley},~\bfnm{Patrick}\binits{P.}}
(\byear{1968}).
\btitle{Convergence of Probability Measures}.
\baddress{New York}: \bpublisher{Wiley}.
\bid{mr={0233396}}
\end{bbook}
%
\endbibitem

\bibitem{Bra07}
%
\begin{bbook}[auto:STB|2011-03-03|12:04:44]
\bauthor{\bsnm{Bradley},~\bfnm{R.~C.}\binits{R.C.}}
(\byear{2007}).
\btitle{Introduction to Strong Mixing Conditions. Vol. I}.
\baddress{Heber City, UT}:
\bpublisher{Kendrick Press}.
\bid{mr={2325294}}
\end{bbook}
%
\endbibitem

\bibitem{BJ94}
%
\begin{barticle}[mr]
\bauthor{\bsnm{Br{\"a}nn{\"a}s},~\bfnm{Kurt}\binits{K.}} \AND
\bauthor{\bsnm{Johansson},~\bfnm{Per}\binits{P.}}
(\byear{1994}).
\btitle{Time series count data regression}.
\bjournal{Comm. Statist. Theory Methods}
\bvolume{23}
\bpages{2907--2925}.
\bid{doi={10.1080/03610929408831424}, issn={0361-0926}, mr={1294013}}
\end{barticle}
%
\endbibitem

\bibitem{CT86}
%
\begin{barticle}[auto:STB|2011-03-03|12:04:44]
\bauthor{\bsnm{Cameron},~\bfnm{A.~C.}\binits{A.C.}} \AND
\bauthor{\bsnm{Trivedi},~\bfnm{P.~K.}\binits{P.K.}}
(\byear{1986}).
\btitle{Econometric models based on count data: Comparison and
applications of
some estimators and tests}.
\bjournal{J. Appl. Econometrics}
\bvolume{1}
\bpages{29--53}.
\end{barticle}
%
\endbibitem

\bibitem{Cox81}
%
\begin{barticle}[mr]
\bauthor{\bsnm{Cox},~\bfnm{D.~R.}\binits{D.R.}}
(\byear{1981}).
\btitle{Statistical analysis of time series: Some recent developments (with discussion)}.
\bjournal{Scand. J. Stat.}
\bvolume{8}
\bpages{93--115}.
\bid{issn={0303-6898}, mr={0623586}}
\bptnote{check related}
\end{barticle}
%
\endbibitem

\bibitem{Cox83}
%
\begin{barticle}[mr]
\bauthor{\bsnm{Cox},~\bfnm{D.~R.}\binits{D.R.}}
(\byear{1983}).
\btitle{Some remarks on overdispersion}.
\bjournal{Biometrika}
\bvolume{70}
\bpages{269--274}.
\bid{doi={10.1093/biomet/70.1.269}, issn={0006-3444}, mr={0742997}}
\end{barticle}
%
\endbibitem

\bibitem{DDS03}
%
\begin{barticle}[mr]
\bauthor{\bsnm{Davis},~\bfnm{Richard~A.}\binits{R.A.}},
\bauthor{\bsnm{Dunsmuir},~\bfnm{William T.~M.}\binits{W.T.M.}} \AND
\bauthor{\bsnm{Streett},~\bfnm{Sarah~B.}\binits{S.B.}}
(\byear{2003}).
\btitle{Observation-driven models for {P}oisson counts}.
\bjournal{Biometrika}
\bvolume{90}
\bpages{777--790}.
\bid{doi={10.1093/biomet/90.4.777}, issn={0006-3444}, mr={2024757}}
\end{barticle}
%
\endbibitem

\bibitem{DDW99}
%
\begin{bincollection}[mr]
\bauthor{\bsnm{Davis},~\bfnm{Richard~A.}\binits{R.A.}},
\bauthor{\bsnm{Wang},~\bfnm{Ying}\binits{Y.}} \AND
\bauthor{\bsnm{Dunsmuir},~\bfnm{William T.~M.}\binits{W.T.M.}}
(\byear{1999}).
\btitle{Modeling time series of count data}.
In \bbooktitle{Asymptotics, Nonparametrics, and Time Series}.
\bseries{Statist. Textbooks Monogr.}
\bvolume{158}
(\beditor{\bfnm{S.}\binits{S.}~\bsnm{Ghosh}}, ed.)
\bpages{63--113}.
\baddress{New York}: \bpublisher{Dekker}.
\bid{mr={1724696}}
\end{bincollection}
%
\endbibitem

\bibitem{DW09}
%
\begin{barticle}[mr]
\bauthor{\bsnm{Davis},~\bfnm{Richard~A.}\binits{R.A.}} \AND
\bauthor{\bsnm{Wu},~\bfnm{Rongning}\binits{R.}}
(\byear{2009}).
\btitle{A negative binomial model for time series of counts}.
\bjournal{Biometrika}
\bvolume{96}
\bpages{735--749}.
\bid{doi={10.1093/biomet/asp029}, issn={0006-3444}, mr={2538769}}
\end{barticle}
%
\endbibitem

\bibitem{DL89}
%
\begin{barticle}[mr]
\bauthor{\bsnm{Dean},~\bfnm{C.}\binits{C.}} \AND
\bauthor{\bsnm{Lawless},~\bfnm{J.~F.}\binits{J.F.}}
(\byear{1989}).
\btitle{Tests for detecting overdispersion in {P}oisson regression models}.
\bjournal{J.~Amer. Statist. Assoc.}
\bvolume{84}
\bpages{467--472}.
\bid{issn={0162-1459}, mr={1010335}}
\end{barticle}
%
\endbibitem

\bibitem{Dea92}
%
\begin{barticle}[auto:STB|2011-03-03|12:04:44]
\bauthor{\bsnm{Dean},~\bfnm{C.~B.}\binits{C.B.}}
(\byear{1992}).
\btitle{Testing for overdispersion in Poisson and binomial regression models}.
\bjournal{J.~Amer. Statist. Assoc.}
\bvolume{87}
\bpages{451--457}.
\end{barticle}
%
\endbibitem

\bibitem{DF99}
%
\begin{barticle}[mr]
\bauthor{\bsnm{Diaconis},~\bfnm{Persi}\binits{P.}} \AND
\bauthor{\bsnm{Freedman},~\bfnm{David}\binits{D.}}
(\byear{1999}).
\btitle{Iterated random functions}.
\bjournal{SIAM Rev.}
\bvolume{41}
\bpages{45--76}.
\bid{doi={10.1137/S0036144598338446}, issn={0036-1445}, mr={1669737}}
\end{barticle}
%
\endbibitem

\bibitem{Dou94}
%
\begin{bbook}[mr]
\bauthor{\bsnm{Doukhan},~\bfnm{Paul}\binits{P.}}
(\byear{1994}).
\btitle{Mixing: Properties and Examples}.
\bseries{Lecture Notes in Statistics}
\bvolume{85}.
\baddress{New York}: \bpublisher{Springer}.
\bid{mr={1312160}}
\end{bbook}
%
\endbibitem

\bibitem{Dur91}
%
\begin{bbook}[mr]
\bauthor{\bsnm{Durrett},~\bfnm{Richard}\binits{R.}}
(\byear{1991}).
\btitle{Probability: Theory and Examples}.
\baddress{Pacific Grove, CA}: \bpublisher{Wadsworth \& Brooks/Cole}.
\bid{mr={1068527}}
\end{bbook}
%
\endbibitem

\bibitem{FLO06}
%
\begin{barticle}[mr]
\bauthor{\bsnm{Ferland},~\bfnm{Ren{\'e}}\binits{R.}},
\bauthor{\bsnm{Latour},~\bfnm{Alain}\binits{A.}} \AND
\bauthor{\bsnm{Oraichi},~\bfnm{Driss}\binits{D.}}
(\byear{2006}).
\btitle{Integer-valued {GARCH} process}.
\bjournal{J. Time Series Anal.}
\bvolume{27}
\bpages{923--942}.
\bid{doi={10.1111/j.1467-9892.2006.00496.x}, issn={0143-9782}, mr={2328548}}
\end{barticle}
%
\endbibitem

\bibitem{FRT09}
%
\begin{barticle}[mr]
\bauthor{\bsnm{Fokianos},~\bfnm{Konstantinos}\binits{K.}},
\bauthor{\bsnm{Rahbek},~\bfnm{Anders}\binits{A.}} \AND
\bauthor{\bsnm{Tj{\o}stheim},~\bfnm{Dag}\binits{D.}}
(\byear{2009}).
\btitle{Poisson autoregression}.
\bjournal{J. Amer. Statist. Assoc.}
\bvolume{104}
\bpages{1430--1439}.
\bid{doi={10.1198/jasa.2009.tm08270}, issn={0162-1459}, mr={2596998}}
\end{barticle}
%
\endbibitem

\bibitem{FT10}
%
\begin{barticle}[auto:STB|2011-03-03|12:04:44]
\bauthor{\bsnm{Fokianos},~\bfnm{K.}\binits{K.}} \AND
\bauthor{\bsnm{Tj{\o}stheim},~\bfnm{D.}\binits{D.}}
(\byear{2011}).
\btitle{Log-linear Poisson
autoregression.}
\bjournal{J. Mult. Anal.} \bvolume{102} \bpages{\mbox{563--578}}.
\end{barticle}
%
\endbibitem

\bibitem{GHTT00}
%
\begin{barticle}[mr]
\bauthor{\bsnm{Grunwald},~\bfnm{Gary~K.}\binits{G.K.}},
\bauthor{\bsnm{Hyndman},~\bfnm{Rob~J.}\binits{R.J.}},
\bauthor{\bsnm{Tedesco},~\bfnm{Leanna}\binits{L.}} \AND
\bauthor{\bsnm{Tweedie},~\bfnm{Richard~L.}\binits{R.L.}}
(\byear{2000}).
\btitle{Non-{G}aussian conditional linear {${\rm AR}(1)$} models}.
\bjournal{Aust. N. Z. J. Stat.}
\bvolume{42}
\bpages{479--495}.
\bid{doi={10.1111/1467-842X.00143}, issn={1369-1473}, mr={1802969}}
\end{barticle}
%
\endbibitem

\bibitem{KF02}
%
\begin{bbook}[mr]
\bauthor{\bsnm{Kedem},~\bfnm{Benjamin}\binits{B.}} \AND
\bauthor{\bsnm{Fokianos},~\bfnm{Konstantinos}\binits{K.}}
(\byear{2002}).
\btitle{Regression Models for Time Series Analysis}.
\baddress{Hoboken, NJ}: \bpublisher{Wiley}.
\bid{doi={10.1002/0471266981}, mr={1933755}}
\end{bbook}
%
\endbibitem

\bibitem{LM89}
%
\begin{barticle}[mr]
\bauthor{\bsnm{Lasota},~\bfnm{Andrzej}\binits{A.}} \AND
\bauthor{\bsnm{Mackey},~\bfnm{Michael~C.}\binits{M.C.}}
(\byear{1989}).
\btitle{Stochastic perturbation of dynamical systems: The weak
convergence of
measures}.
\bjournal{J. Math. Anal. Appl.}
\bvolume{138}
\bpages{232--248}.
\bid{doi={10.1016/0022-247X(89)90333-8}, issn={0022-247X}, mr={0988333}}
\end{barticle}
%
\endbibitem

\bibitem{Lee86}
%
\begin{barticle}[mr]
\bauthor{\bsnm{Lee},~\bfnm{Lung~Fei}\binits{L.F.}}
(\byear{1986}).
\btitle{Specification test for {P}oisson regression models}.
\bjournal{Internat. Econom. Rev.}
\bvolume{27}
\bpages{689--706}.
\bid{doi={10.2307/2526689}, issn={0020-6598}, mr={0863069}}
\end{barticle}
%
\endbibitem

\bibitem{MT93}
%
\begin{bbook}[mr]
\bauthor{\bsnm{Meyn},~\bfnm{S.~P.}\binits{S.P.}} \AND
\bauthor{\bsnm{Tweedie},~\bfnm{R.~L.}\binits{R.L.}}
(\byear{1993}).
\btitle{Markov Chains and Stochastic Stability}.
\baddress{London}: \bpublisher{Springer.}
\bid{mr={1287609}}
\end{bbook}
%
\endbibitem

\bibitem{Pol84}
%
\begin{bbook}[mr]
\bauthor{\bsnm{Pollard},~\bfnm{David}\binits{D.}}
(\byear{1984}).
\btitle{Convergence of Stochastic Processes}.
\baddress{New York}: \bpublisher{Springer}.
\bid{mr={0762984}}
\end{bbook}
%
\endbibitem

\bibitem{Ros80}
%
\begin{barticle}[mr]
\bauthor{\bsnm{Rosenblatt},~\bfnm{M.}\binits{M.}}
(\byear{1980}).
\btitle{Linear processes and bispectra}.
\bjournal{J. Appl. Probab.}
\bvolume{17}
\bpages{265--270}.
\bid{issn={0021-9002}, mr={0557456}}
\end{barticle}
%
\endbibitem

\bibitem{RS00}
%
\begin{bmisc}[auto:STB|2011-03-03|12:04:44]
\bauthor{\bsnm{Rydberg},~\bfnm{T.~H.}\binits{T.H.}} \AND
\bauthor{\bsnm{Shephard},~\bfnm{N.}\binits{N.}}
(\byear{2000}).
\bhowpublished{BIN models for trade-by-trade data. Modelling the number
of trades in a fixed interval of time. In \textit{World Conference Econometric Society,
2000, Seattle}. Contributed Paper 0740}.
\end{bmisc}
%
\endbibitem

\bibitem{Str00}
%
\begin{bmisc}[mr]
\bauthor{\bsnm{Streett},~\bfnm{Sarah~Burns}\binits{S.B.}}
(\byear{2000}).
\bhowpublished{Some observation driven models for time series.
Ph.D. thesis, Dept. Statistics, Colorado State Univ}.
\bid{mr={2701787}}
\end{bmisc}
%
\endbibitem

\bibitem{ZB08}
%
\begin{barticle}[mr]
\bauthor{\bsnm{Zheng},~\bfnm{Haitao}\binits{H.}} \AND
\bauthor{\bsnm{Basawa},~\bfnm{Ishwar~V.}\binits{I.V.}}
(\byear{2008}).
\btitle{First-order observation-driven integer-valued autoregressive
processes}.
\bjournal{Statist. Probab. Lett.}
\bvolume{78}
\bpages{1--9}.
\bid{doi={10.1016/j.spl.2007.04.017}, issn={0167-7152}, mr={2381267}}
\end{barticle}
%
\endbibitem

\end{thebibliography}
\end{document}